\documentclass[12pt,letterpaper]{article}
\usepackage[hmargin=1in,vmargin=1in]{geometry}

\usepackage{amsmath}
\usepackage{amsfonts}
\usepackage{amsthm}
\usepackage{amssymb}
\usepackage{array}
\usepackage{caption}
\usepackage{color}
\usepackage{graphics}
\usepackage{graphicx}
\usepackage{hyperref}
\usepackage{verbatim}

\newtheorem*{tp}{Tanya's Program}

\begin{document}

\title{Hat Puzzles}
\author{Tanya Khovanova}
\maketitle

\begin{abstract}
This paper serves as the announcement of my program---a joke version of the Langlands Program. In connection with this program, I discuss an old hat puzzle, introduce a new hat puzzle, and offer a puzzle for the reader.
\end{abstract}

\section{Tanya's Program}
\label{sec:tp}

This is one of my favorite jokes:

\begin{quote}
Three logicians walk into a bar. The waitress asks, ``Do you all want beer?''\\
The first logician answers, ``I do not know.''\\
The second logician answers, ``I do not know.''\\
The third logician answers, ``Yes.''
\end{quote}

This joke reminds me of hat puzzles. In the joke each logician knows whether or not s/he wants a beer, but doesn't know what the others want to drink. In hat puzzles logicians know the colors of the hats on others' heads, but not the color of their own hats.

Here is a hat puzzle where the logicians provide the same answers as in the beer joke. Three logicians wearing hats walk into a bar. They know that the hats were placed on their heads from the set of hats in Figure~\ref{fig:hatset}. The total number of available red hats was three, and the total number of available blue hats was two.

\begin{center}
\begin{tabular}{@{}ccccc@{}}
\includegraphics[height=9mm]{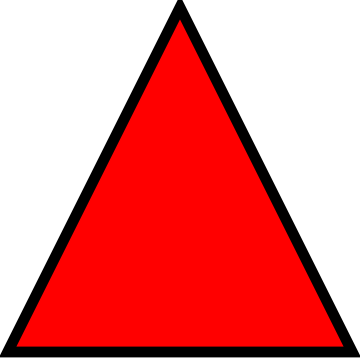} &
\includegraphics[height=9mm]{RHat.png} &
\includegraphics[height=9mm]{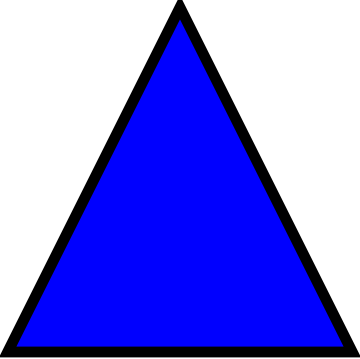} &
\includegraphics[height=9mm]{BHat.png} &
\includegraphics[height=9mm]{RHat.png}\\
\end{tabular}
\captionof{figure}{The set of available hats}\label{fig:hatset}
\end{center}

\begin{quote}
Three logicians walk into a bar. The waitress asks, ``Do you know the color of your own hat?''\\
The first logician answers, ``I do not know.''\\
The second logician answers, ``I do not know.''\\
The third logician answers, ``Yes.''
\end{quote}

The question is: What is the color of the third logician's hat?

This process of converting jokes to puzzles reminds me of the Langlands Program, which tries to unite different parts of mathematics. I would like to unite jokes and puzzles. So here I announce my own program:

\begin{tp}
Find a way to convert jokes into puzzles and puzzles into jokes.
\end{tp}

\section{Old Hat Puzzle}

Enough jokes. Let us move on to hat puzzles. The following is one of my favorite hat puzzles:

\begin{quote}
Some logicians are put to the test. They will have time to develop a strategy before the test.

Here is the test: The logicians stand in a line, one behind the other (as in Figure~\ref{fig:oldpuzzle}), so that the last person in line can see everyone else. Each wears a red or a blue hat and, as is usual in such puzzles, none of them knows the color of his or her own hat. Also, unlike in the puzzle above, there is no limit on the number of available hats of each color. The logicians can only see the colors of the hats on all of the people in front of them. Then, one at a time, according to the strategy they have agreed in advance, each logician guesses out loud the color of the hat on his or her own head. They are not permitted to give any other information: they can't even sneeze, fart or poke the person in front of them. They can only convey one bit of information: red or blue.

Their goal is to maximize the number of logicians who are guaranteed to guess correctly.
\end{quote}

\begin{center}
\begin{tabular}{@{}cccccc@{}}
\includegraphics[height=15mm]{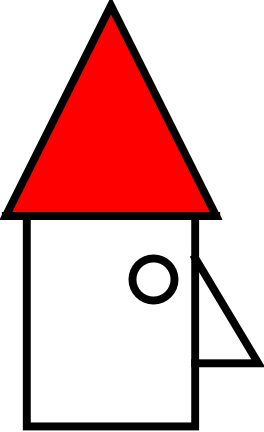} &
\includegraphics[height=15mm]{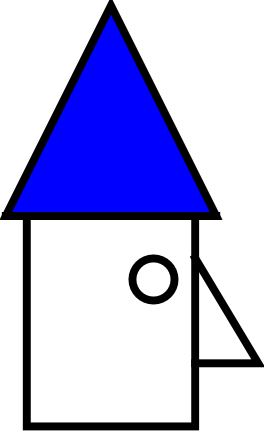}  &
\includegraphics[height=15mm]{RedHat.png} &
\includegraphics[height=15mm]{BlueHat.png}&
\includegraphics[height=15mm]{RedHat.png} &
\includegraphics[height=15mm]{RedHat.png}\\ 
\end{tabular}
\captionof{figure}{An old hat puzzle}\label{fig:oldpuzzle}
\end{center}

I don't know the origins of this puzzle. It appeared in the 23rd All-Russian Mathematical Olympiad in 1997 as folklore. You can also find it in Peter Winkler's math puzzle book for gourmands \cite{PW} and in an article by Ezra Brown and James Tanton with a dozen hat puzzles \cite{BT}.

\subsection{Two Logicians}

Let us start with two logicians in line as in Figure~\ref{fig:two}. 

\begin{center}
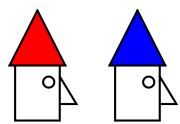

\begin{tabular}{@{}cc@{}}
\includegraphics[height=15mm]{RedHat.png} &
\includegraphics[height=15mm]{BlueHat.png} \\ 
\end{tabular}
\captionof{figure}{Two logicians}\label{fig:two}
\end{center}

I teach a Math Competitions class at the Advanced Math and Sciences Academy Charter School. I gave my students this puzzle. In class we started with two logicians. It is impossible to guarantee that the last person can correctly guess the color of his or her hat, because no one sees it. The best we can hope to do is to guarantee that one person guesses the color correctly, namely the first logician in line. Most of my students solve this case very quickly. They realize that the last person in line should name the color of the hat in front of him/her, and the front logician should repeat this color.

\subsection{Three Logicians}

Three logicians (See Figure~\ref{fig:three}) give my students a moment's pause.

\begin{center}
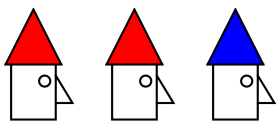

\begin{tabular}{@{}ccc@{}}
\includegraphics[height=15mm]{RedHat.png} &
\includegraphics[height=15mm]{RedHat.png} &
\includegraphics[height=15mm]{BlueHat.png} \\ 
\end{tabular}
\captionof{figure}{Three logicians}\label{fig:three}
\end{center}

The students start by suggesting that the last person names either the color of the hat in front of him/her, or the color of the first person's hat. They soon realize that these strategies do not work. The last person needs to combine together all the colors s/he sees. Then they discover the solution to guarantee that all but the last logician guess correctly: the last person signals whether the two people in front have the same color or different color hats.

For example, the logicians can agree beforehand that ``red'' means the same color and ``blue'' means different colors. In the situation depicted in Figure~\ref{fig:three}, the last person sees different colors, so s/he says blue. The second to last logician now knows that the colors are different, and s/he sees blue, so she says ``red.'' Now the first logician in line knows that his/her hat color is different from the middle person's color. Since the middle person said ``red,'' the first person says ``blue.''

\subsection{Many Logicians}

The case with four logicians takes some more thinking. The students usually solve the four-logician case together with the any-number-of-logicians case (See Figure~\ref{fig:many}). Surprisingly, not more than one logician needs to be mistaken.

\begin{center}
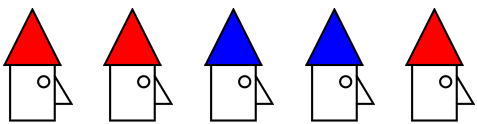

\begin{tabular}{@{}ccccc@{}}
\includegraphics[height=15mm]{RedHat.png} &
\includegraphics[height=15mm]{RedHat.png} &
\includegraphics[height=15mm]{BlueHat.png} &
\includegraphics[height=15mm]{BlueHat.png} &
\includegraphics[height=15mm]{RedHat.png}\\ 
\end{tabular}
\captionof{figure}{Many logicians}\label{fig:many}
\end{center}

By this time the students understand that the last logician must start. S/he is the only one who can't have any information about his/her hat's color. So the last logician is doomed and uses this as an opportunity to pass information to everyone else. It is also clear that the logicians should state their colors from back to front. This way all logicians---except for the last one---are in the same situation. By the time the logicians have had to guess their own color, they know all the other colors: they see all the colors in front of them and have heard all the colors of the logicians behind them. 

The solution: The last logician has to announce the parity of the number of red hats. For example, the last person says ``Red'' if the number of red hats s/he sees is even. After that the other logicians sum up the number of times they see or hear the word ``Red,'' and say ``Red'' if this number is even and ``Blue'' if it is not.

I make my students re-enact this test. If more than one student is wrong, I chop off all of their heads (See Figure~\ref{fig:off}). Just kidding.

\begin{center}
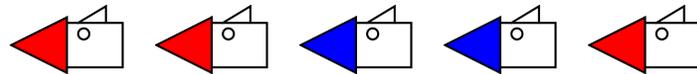

\begin{tabular}{@{}ccccc@{}}
\includegraphics[height=15mm, angle=90]{RedHat.png} &
\includegraphics[height=15mm, angle=90]{RedHat.png} &
\includegraphics[height=15mm, angle=90]{BlueHat.png} &
\includegraphics[height=15mm, angle=90]{BlueHat.png} &
\includegraphics[height=15mm, angle=90]{RedHat.png}\\ 
\end{tabular}
\captionof{figure}{Heads chopped off}\label{fig:off}
\end{center}

\subsection{Many Colors}

It is natural to generalize this problem to many colors. The thee-color version was suggested by Konstantin Knop for the 23rd All-Russian Mathematical Olympiad in 1997. You can also find the 100-color version in the above-mentioned paper with a dozen hat puzzles \cite{BT}.

We have the same logicians standing in a line and now they might have a hat of any color---not just red or blue (See Figure~\ref{fig:manycolors}). The important thing is that the set of hat colors is finite and needs to be known in advance.

\begin{center}
\begin{tabular}{@{}cccccc@{}}
\includegraphics[height=15mm]{RedHat.png} &
\includegraphics[height=15mm]{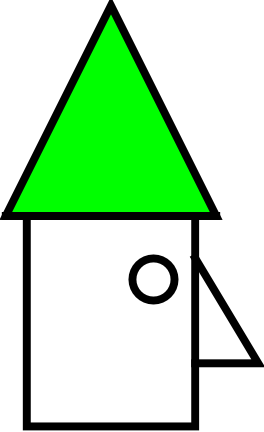} &
\includegraphics[height=15mm]{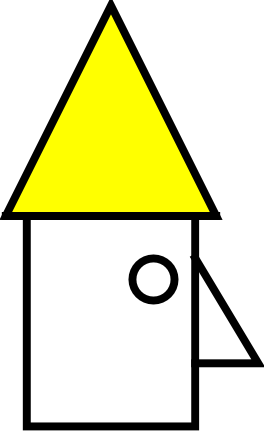} &
\includegraphics[height=15mm]{BlueHat.png} &
\includegraphics[height=15mm]{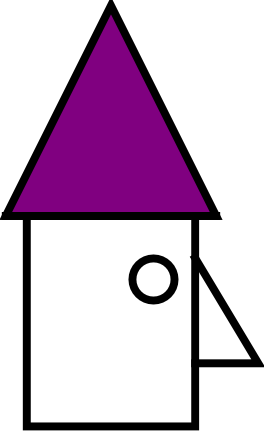} &
\includegraphics[height=15mm]{GreenHat.png}\\ 
\end{tabular}
\captionof{figure}{Many colors}\label{fig:manycolors}
\end{center}

It might sound surprising, but again the logicians can guarantee that not more than one person, namely the last person, is mistaken. The strategy is similar. Logicians replace colors by numbers from 0 to $N-1$, where $N$ is the number of possible colors. The last person sums the colors modulo $N$ and announces his/her own color so that the total sum is zero. As before, the other logicians state their colors from back to front. Each logician sums up all the numbers, that is colors, that s/he sees and hears. Then s/he calculates his/her own color by assuming that the total sum modulo $N$ is zero.

\section{New Hat Puzzle}

This new variation is very recent and very beautiful. It was invented by Konstantin Knop and Alexander Shapovalov and appeared (with different wording) in March of 2013 at the Tournament of the Towns:

\begin{quote}
The same logicians stand in line, one behind the other, so that the last person in line sees everyone else. They were previously shown all the hats from which the hats on their heads have been chosen. Every available hat is a different color and there is one more hat than the number of logicians (see Figure~\ref{fig:distinctcolors}). As before, in any order they agreed earlier, each logician guesses out loud the color of the hat on his/her own head, but is not permitted to signal anything else. No one can repeat a color that has already been announced. The logicians have time to design the strategy before the test and they need to maximize the number of people who are guaranteed to guess correctly. What should that strategy be?
\end{quote}

\begin{center}
\begin{tabular}{@{}cccccc@{}}
\includegraphics[height=7.5mm]{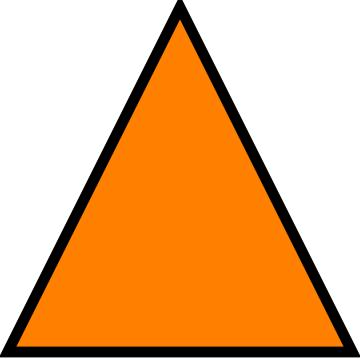} &
\includegraphics[height=15mm]{YellowHat.png} &
\includegraphics[height=15mm]{RedHat.png} &
\includegraphics[height=15mm]{PurpleHat.png} &
\includegraphics[height=15mm]{BlueHat.png} &
\includegraphics[height=15mm]{GreenHat.png}\\ 
\end{tabular}
\captionof{figure}{Distinct colors}\label{fig:distinctcolors}
\end{center}

It is natural to try to reuse the solution we know for the previous puzzle. Let us say that the number of hats is $N$ and the number of logicians is $N-1$. Suppose the last person sums the colors modulo $N$ and announces his/her own color so that the total sum is zero. The problem is that s/he might announce the color of one of the logicians in front of him/her. 

Suppose the last logician says ``Red.'' When it's time for the logician who is actually wearing a red hat to speak, s/he will realize that his/her hat is red, and oops, s/he can't repeat this color. If s/he says a different color randomly, s/he will throw off the strategy and mislead everyone in front.

Is there a way to rescue this strategy so that only a small number of logicians might be mistaken? It appears that there is. Suppose the last logician says ``Red,'' and the logician in the red hat says ``Blue.'' The blue color wasn't announced yet so the logician in the blue hat is most probably way in front of the red-hatted logician. That means that all the logicians between the red-hatted and the blue-hatted logicians will know that something wrong happened. They will not be misled. They will realize that the logician who said ``Blue'' said so because s/he were not allowed to say his/her own color. So that color must be red. That means all the logicians in between can still figure out their colors.

Here is the strategy where not more than three people will name a wrong color:

The last person announces the color so that his/her color and the sum of the colors s/he sees modulo $N$ is zero. If there is a logician with this hat color, s/he says the hat color of the first person. This way the three people who might be mistaken are the last one, the logician with the hat color announced by the last person and the first person in line. Everyone else is guaranteed to be correct.

This is quite a good solution, but there is a completely different solution that guarantees that not more than \textbf{one} person is mistaken. I invite you to try and solve it yourself. If you can't do it, you can find the solution  in my paper devoted to this puzzle \cite{TK}.

\section{A Puzzle for the Reader}

Convert the above puzzles into jokes!


\begin{thebibliography}{9}

\bibitem{BT} Ezra Brown and James Tanton, A Dozen Hat Problems, Math Horizons, 16 No. 8 (2009) 22--25.

\bibitem{TK} Tanya Khovanova, A Line of Sages, to appear in The Mathematical Intelligencer, 3, 2014.

\bibitem{PW} Peter Winkler, \emph{Mathematical Puzzles: A Connoisseur's Collection}, A.K. Peters, Natick, MA, 2004.



\end{thebibliography}
\end{document}